\documentclass{elsarticle}
\usepackage{amsmath,amssymb}
\newcommand{\E}{{\bf E}}
\renewcommand{\P}{{\bf P}}
\newcommand{\sgn}{{\rm sgn}}

\begin{document}

\title{On Cox-Kemperman moment inequalities for independent centered random variables}

\author{P.S.Ruzankin\footnote{Sobolev Institute of Mathematics, Novosibirsk State University, pr. Ak. Koptyuga, 4, Novosibirsk, Russia,
ruzankin@math.nsc.ru}}

\date{}
%



\begin{abstract}
In 1983 Cox and Kemperman proved that
$\E f(\xi)+ \E f(\eta) \le \E f(\xi+\eta)$ for all functions $f$, such that $f(0)=0$ and the second derivative
$f''(y)$ is convex, and all independent centered random variables $\xi$ and $\eta$ satisfying
certain moment restrictions.
We show that the minimal moment restrictions are sufficient for the inequality to be valid,
and write out
a less restrictive condition on $f$ for the inequality to hold. 
%

Besides, Cox and Kemperman (1983)
found out the optimal constants $A_\rho$ and $B_\rho$ for the
inequalities
$A_\rho (\E |\xi|^\rho + \E |\eta|^\rho) \le
\E |\xi + \eta|^\rho \le
B_\rho (\E |\xi|^\rho + \E |\eta|^\rho)
$, where $\rho\ge1$, $\xi$ and $\eta$ are independent centered random variables. We write out similar sharp inequalities for symmetric random variables.


\end{abstract}

\maketitle

\section{Introduction and formulation of the results}

Cox and Kemperman have proved the following theorem:

{\bf Theorem A} [Cox and Kemperman, 1983]. {\it

Let random variables $\xi$ and $\eta$ be such that
\begin{equation}\label{ck0}
\E (\xi|\eta)=0,\ \ \E(\eta|\xi)=0 \ \mbox{ a.s.}
\end{equation}

Then, for each $\rho \ge 1$, the following inequalities hold:
\begin{equation}\label{ck1}
2^{\rho-2} (\E |\xi|^\rho + \E |\eta|^\rho) \le
\E |\xi + \eta|^\rho \le
\left(\max_{0\le z\le 1} \psi(\rho,z)\right) (\E |\xi|^\rho + \E |\eta|^\rho)
\quad\mbox{ if } 1\le \rho\le 2,
\end{equation}
\begin{equation}\label{ck2}
\left(\min_{0\le z\le 1} \psi(\rho,z)\right) (\E |\xi|^\rho + \E |\eta|^\rho) \le
\E |\xi + \eta|^\rho \le
2^{\rho-2} (\E |\xi|^\rho + \E |\eta|^\rho)
\quad\mbox{ if } 2\le \rho\le 3,
\end{equation}
\begin{equation}\label{ck3}
\E |\xi|^\rho + \E |\eta|^\rho \le
\E |\xi + \eta|^\rho \le
2^{\rho-2} (\E |\xi|^\rho + \E |\eta|^\rho)
\quad\mbox{ if } \rho\ge 3
\end{equation}
whenever $\E |\xi|^\rho<\infty$ and $\E |\eta|^\rho<\infty$,
where
$$\psi(\rho,z)= 2^{\rho-1}(z+z^{\rho-1}+(1-z)^\rho \big)/
\big( (1+z)(1+z^{\rho-1}) .$$

All the estimates in $(\ref{ck1})-(\ref{ck3})$ for $\E |\xi + \eta|^\rho$ are sharp in the sense that, for each
inequality, there exist distributions of $\xi$ and $\eta$, such that
 $\xi\ne0$
and the inequality turns into equality for independent $\xi$ and $\eta$ with these distributions.
}

This theorem does not consider the case $0<\rho<1$ because in the case the sharp inequalities are trivial ones: $0 \le
\E |\xi + \eta|^\rho \le
\E |\xi|^\rho + \E |\eta|^\rho$.

Besides, Cox and Kemperman (1983) noted that, for i.i.d. $\xi$ and $\eta$
having a symmetric two-point distribution,
\begin{equation}\label{ck9}
\E |\xi + \eta|^\rho=2^{\rho-2}(\E |\xi|^\rho + \E |\eta|^\rho )
\quad\mbox{ for all } \rho>0.
\end{equation}

It is also known (e.g. see Rosenthal (1972)) that,
for symmetric independent random variables $\xi$ and $\eta$,
\begin{equation}\label{rt}
\E |\xi|^\rho + \E |\eta|^\rho \le
\E |\xi + \eta|^\rho
\quad\mbox{ if } \rho\ge 2
\end{equation}
when $\E |\xi|^\rho<\infty$ and $\E |\eta|^\rho<\infty$.

Besides of estimates for expectations of power functions, there have been obtained certain inequalities for expectations for some other classes of functions.

{\bf Theorem B} [Cox and Kemperman, 1983]. {\it
Let a function $f$ on the real line be such that $f(0)\le 0$ and
the second derivative $f''(y)$ exists for all $y$ and is convex.
Let random variables $\xi$ and $\eta$ satisfy condition
$(\ref{ck0})$ and be such that
\begin{equation}\label{ock}
\E |f'(\xi)|<\infty,\quad\E|f'(\eta)|<\infty.
\end{equation}

Then
$$\E f(\xi)+ \E f(\eta) \le \E f(\xi+\eta) .$$
}

Note that by Theorem~E below condition (\ref{ock}) can be omitted for independent $\xi$ and $\eta$.

Appendix A below contains a simple proof of this theorem proposed by Borisov. (The proof is valid under certain moment restrictions which can also be omitted for independent $\xi$ and $\eta$ by Theorem~E.)

Note that the function $f(y)=|y|^\rho$ satisfies the conditions of this theorem only if
$\rho=2$ or $\rho\ge 3$.

Utev  has obtained the following corresponding result:

{\bf Theorem C} [Utev, 1985]. {\it
Let a function $f$ on the real line have $f''(y)$ for all $y$. Then the following three
conditions are equivalent

$1)$. $f''$ is convex.

$2)$. For all independent symmetric random variables $\xi$ and $\eta$ and for all $y$,
the inequality
$$ \E f(y+\xi) + \E f(y+\eta) \le  f(y)+\E f(y+\xi+\eta)$$
holds whenever these expectations exist.

$3)$. For all independent bounded random variables $\xi$ and $\eta$, such that
$\E \xi=\E \eta=0$, and for all~$y$,
the inequality
$$ \E f(y+\xi) + \E f(y+\eta)\le f(y)+\E f(y+\xi+\eta)$$
holds.
}

In fact, Utev formulated and proved this theorem for Hilbert space - valued random variables. However such spaces will not be discussed in the present paper.

Note that the restriction that the random variables $\xi$ and $\eta$ in condition 3)
be bounded can be omitted by Theorem~E below.

Another class of functions is considered in the following statement.

{\bf Theorem D} [Figiel, Hitczenko, Johnson, Schechtman, and Zinn, 1997].
{\it
Let a function $f$, $f(0)\le0$, be even and such that the function $y\mapsto f(\sqrt{|y|})$
is convex.
Then
$$\E f(\xi)+ \E f(\eta) \le \E f(\xi+\eta) $$
for all random variables $\xi$ and $\eta$, such that the conditional distribution of
$\eta$ under the condition $\xi=x$ is symmetric for all $x$.
}

Note that the function $|y|^\rho$ satisfies conditions of this theorem only if
$\rho\ge 2$.

We will call a function $f$ on the real line \emph{twice differentiable} if
$f'(y)$ exists for all $y$, $f''(y)$ exists for almost all (with respect to the Lebesgue measure)
$y$, and, for all $a<b$, $f'(b)-f'(a)=\int_a^b f''(y) dy$.


{\bf Theorem 1.} {\it
Let a function $f$ be twice differentiable, $f(0)\le0$, and the function $f''(t)+f''(-t)$ be nondecreasing for
$t>0$.

Then, for all independent symmetric random variables  $\xi$ and $\eta$,
$$\E f(\xi)+ \E f(\eta) \le \E f(\xi+\eta)$$
whenever the expectations exist.
}

A useful corollary of the theorem is the following one.

{\bf Corollary 1.} {\it
For independent symmetric random variables $\xi$ and $\eta$,
the following inequalities are valid:
\begin{equation}\label{c11}
2^{\rho-2} (\E |\xi|^\rho + \E |\eta|^\rho) \le
\E |\xi + \eta|^\rho \le
\E |\xi|^\rho + \E |\eta|^\rho
\quad\mbox{ if } 0< \rho\le 2,
\end{equation}
\begin{equation}\label{c12}
\E |\xi|^\rho + \E |\eta|^\rho \le
\E |\xi + \eta|^\rho \le 2^{\rho-2} (
\E |\xi|^\rho + \E |\eta|^\rho)
\quad\mbox{ if } \rho\ge 2
\end{equation}
when $\E |\xi|^\rho<\infty$ and $\E |\eta|^\rho<\infty$.
The four estimates for $\E |\xi + \eta|^\rho$ are sharp in the sense that, for each inequality, there exist distributions
of $\xi$ and $\eta$, such that $\xi\ne 0$ and the inequality turns into equality for these distributions.
}

Note that (\ref{c12}) and (\ref{c11}) for $\rho=1$ follow directly from Theorem A and relations (\ref{ck9}) and (\ref{rt}). Besides, the left inequality in (\ref{c11}) for $1\le\rho\le 2$ follows directly from Theorem A and relation (\ref{ck9}).

Note also that the right inequality in (\ref{c11}) for $0<\rho\le 1$ is trivial because
$|\alpha+\beta|^\rho\le|\alpha|^\rho+|\beta|^\rho$ for any real numbers $\alpha$ and $\beta$,
$0<\rho\le 1$.


{\bf Theorem 2.}
{\it
Let a function $f$ be twice differentiable and such that $f(0)\le 0$ and
\begin{equation}\label{t30}
f''(-\alpha)+f''(\beta)\ge f''(-\alpha+\gamma)+
f''(\beta-\gamma),
\end{equation}
for any $\alpha>0$, $\beta>0$, $0<\gamma<\alpha+\beta$ such that
$f''$ is defined at the points $-\alpha$, $\beta$, $-\alpha+\gamma$ and
$\beta-\gamma$.

Then, for all independent centered random variables  $\xi$ and $\eta$,
$$\E f(\xi)+ \E f(\eta) \le \E f(\xi+\eta)$$
whenever the expectations exist.
}

Note that if $f''$ is convex then it satisfies the condition (\ref{t30}).
But the class of functions subject to condition (\ref{t30})
is wider than the class of functions with convex second derivative.

For instance, if $f''(y)=\lfloor h(y)\rfloor$, where $h(y)$ is nonnegative and convex, $h(0)=0$,
then $f(y)$ satisfies condition (\ref{t30}). $\lfloor \cdot \rfloor$ denotes integer part of a number.

As another example, we can take $f''(y)=-y$ if $y<1$, $f''(y)= \lfloor y \rfloor -
(y-\lfloor y \rfloor)$ if $y\ge 1$. Such $f(y)$ also satisfies (\ref{t30}).

{\bf Remark.}
In Theorem 2, if $\xi+\eta \in [-B,C]$ a.s. then it is sufficient to require that the function $f$ satisfy condition (\ref{t30}) only for $\alpha$, $\beta$ lying in $(-B,C)$.

In Theorem 1, if $\xi+\eta\le C$ a.s. then it is sufficient
to require that
$f''(t)+f''(-t)$ be nondecreasing for
$0<t<C$.

The proof of Theorem~1 is based on the fact that any symmetric distribution can be ``decomposed'' into a mixture of symmetric  distributions (e.g. see Figiel, Hitczenko et al. (1997)).
As for Theorem~2, any centered distribution can be ``decomposed'' into a mixture of two-point centered distributions (e.g. see Pinelis (2009) and references therein). So one can prove the corresponding inequalities for two-point distributions only:

{\bf Theorem E.}
{\it
Let a function $g$ of \ $m+n$ arguments, $m\ge 0 $ and $n\ge 0$, be such that
$$
\E g(\xi_1,...,\xi_m,\eta_1,...,\eta_n) \ge 0
$$
for all independent random variables $\xi_j$ and $\eta_j$, where
each of the random variables $\xi_j$ has a centered two-point distribution or equals zero,
and each of $\eta_j$ has a symmetric two-point distribution or equals zero.

Then, for all independent random variables $\xi_1,...,\xi_m$, $\eta_1,...,\eta_n$, where
the random variables $\xi_j$ are centered
and  $\eta_j$ are symmetric, the following inequality is valid:
$$
\E g(\xi_1,...,\xi_m,\eta_1,...,\eta_n) \ge 0
$$
whenever the expectation exists.
}

%
%
%
%

\section{Proofs}

\subsection{Proof of Theorem E}

For the sake of convenience we give here the proof of Theorem~E, but for the case $m=2$, $n=0$ only.
The case of arbitrary $m$ and $n$ can be considered analogously.


Denote $\xi=\xi_1$, $\eta=\xi_2$. If $\xi$ and $\eta$ have centered two-point distributions,
$\xi$ takes values $-a,b$ and $\eta$ takes values $-c,d$ then
$$\P(\xi=-a)=b/(a+b),\ \ \P(\xi=b)=a/(a+b),
\ \ \P(\eta=-c)=d/(c+d),\ \ \P(\xi=d)=c/(c+d).
$$
Thus we have
\begin{equation}\label{ef}
\E f(\xi,\eta)=\frac{1}{(a+b)(c+d)} \big(
bd\, f(-a,-c)+
bc\, f(-a,d)+
ad\, f(b,-c)+
ac\, f(b,d)   \big)\ge 0
\end{equation}
for all $a,b,c,d>0$.

Now let $\xi$ and $\eta$ have arbitrary centered distributions such that
$\P(\xi\ne 0)=\P(\eta\ne 0)=1$. Put
$$p=\P (\xi>0),\quad
F_{\xi}(u)=\P (\xi < u) - (1-p), \quad G_{\xi}(u)=\P (-\xi <u) - p.$$
Put also
$$s(y)=\int_0^y F_\xi^{(-1)}(u) du,\quad
t(x)=\int_0^x G_\xi^{(-1)}(u) du,$$
where $F_{\xi}^{(-1)}(u):=\sup\{x: F_{\xi}(x)<u\}$ is the quantile transformation of $F_{\xi}$,

Then $s(y)$ and $t(x)$ are (strictly) increasing continuous functions
on $[0,p]$ and $[0,1-p]$, respectively, and
$$s(p)=\E \max\{0,\xi\}=\E \max\{0,-\xi\}=t(1-p).
$$
Put
$$z(y)=t^{-1}(s(y)).$$
We have $t(z(y))=s(y)$, hence $dt(z(y))=ds(y)$ which can be rewritten as
$$G_\xi^{(-1)}(z(y))\,dz(y)= F_\xi^{(-1)}(y)\,dy .$$
For a function $h$,
$$\E h(\xi)= \int_0^p h\big(F_\xi^{(-1)}(y)\big)\,dy +
\int_0^{1-p} h\big(-G_\xi^{(-1)}(x)\big)\,dx.
$$
Substituting $x=z(y)$ into the last integral yields
\begin{equation}\label{aae}
\E h(\xi)= \int_0^p \left( h\big(F_\xi^{(-1)}(y)\big)
+h\big(-G_\xi^{(-1)}(z(y)\big)\frac{F_\xi^{(-1)}(y)}{G_\xi^{(-1)}(z(y))}
\right)
\,dy.
\end{equation}

Let us introduce the same notations for $\eta$. Put
$$q=\P (\eta>0),
\quad
F_{\eta}(t)=\P (\eta < u) - (1-q), \quad G_{\eta}(t)=\P (-\eta <u) - q,$$
and let $w(v)$ be defined by the relations
$$w(0)=0,\quad \quad G_{\eta}^{(-1)}(w(v))dw(v)=F_{\eta}^{(-1)}(v)dv.
$$

Using the above notations we can write
$$\E g(\xi,\eta)=\int_0^q \int_0^p \psi(y,v) dy\, dv,$$
where
$$
\psi(y,v):= g\big(F_{\xi}^{(-1)}(y),F_{\eta}^{(-1)}(v)\big) +
g\big(F_{\xi}^{(-1)}(y),-G_{\eta}^{(-1)}(w(v))\big) \frac{F_{\eta}^{(-1)}(v)}{G_{\eta}^{(-1)}(w(v))}+
$$
$$
g\big(-G_{\xi}^{(-1)}(z(y)),F_{\eta}^{(-1)}(v)\big) \frac{F_{\xi}^{(-1)}(y)}{G_{\xi}^{(-1)}(z(y))}
+
g\big(-G_{\xi}^{(-1)}(z(y)),-G_{\eta}^{(-1)}(w(v))\big) \frac{F_{\xi}^{(-1)}(y)}{G_{\xi}^{(-1)}(z(y))}
\frac{F_{\eta}^{(-1)}(v)}{G_{\eta}^{(-1)}(w(v))},
$$
and $\psi(y,v)\ge 0 $ by relation (\ref{ef}).

We have proved the statement of the theorem for the case $\P(\xi\ne 0)=\P(\eta\ne 0)=1$.
The case $\P(\xi= 0)>0$ or $\P(\eta= 0)>0$ can be easily dealt with using mixtures of zero
and nonzero distributions.

\subsection{Proof of Theorem 2}

Without loss of generality we can assume $f(0)=0$.

By Theorem E it is sufficient to prove the statement of Theorem~2 for all $\xi$ and $\eta$ with centered two-point distributions.

Take $\xi\in\{-a,b\}$, $\eta\in\{-c,d\}$, where $a,b,c,d>0$.
We have
$$
{\bf E}f (\xi+\eta) - \E f(\xi) - \E f(\eta) = \frac{1}{(a+b)(c+d)} \phi(a,b,c,d),$$
where
$$
\begin{array}{rcrl}
\phi(r,s,t,u)&=&
su& \big( f(-r-t)-f(-r)-f(-t) \big)\\
&+&st& \big( f(-r+u)-f(-r)-f(u)  \big)\\
&+&ru& \big( f(s-t) -f(s)-f(-t)  \big)\\
&+&rt& \big( f(s+u) -f(s)-f(u)   \big).
\end{array}
$$
Note that $\phi(r,s,t,u)=0$ if $r=0$ or $s=0$ or $t=0$ or $u=0$.
Moreover,
$$
\frac{\partial^4}{\partial r \partial s \partial t \partial u}\phi(r,s,t,u)=
f''(-r-t)+f''(s+u)-f''(-r+u)-f''(s-t)\ge 0
$$
for positive $r,s,t,u$ by condition (\ref{t30}).

Therefore
$$0\le \int_0^d \int_0^c \int_0^b \int_0^a
\frac{\partial^4}{\partial r \partial s \partial t \partial u}\phi(r,s,t,u)
\,dr\, ds\, dt\, du =
$$
$$
\sum_{ r\in\{0,a\},s\in\{0,b\},t\in\{0,c\},u\in\{0,d\} }
(-1)^{\sgn r+\sgn s+\sgn t+\sgn u} \phi(r,s,t,u)=
\phi(a,b,c,d),
$$
where $\sgn r=0$ if $r=0$ and $\sgn r=1 $ if $r>0$. Thus
$$\phi(a,b,c,d)\ge 0,$$
and hence the theorem is proved.

\subsection{Proof of Theorem 1}

The proof is analogous to that of Theorem 2.

Without loss of generality we can assume $f(0)=0$.

By Theorem E it is sufficient to prove the statement of Theorem~1 for all $\xi$ and $\eta$ with symmetric two-point distributions.

Take $\xi\in\{-a,a\}$, $\eta\in\{-b,b\}$, where $a,b>0$.
We have
$$
{\bf E}f (\xi+\eta) - \E f(\xi) - \E f(\eta) = \frac{1}{4} \phi(a,b),$$
where
$$
\phi(r,s)=
 \big( f(-r-s)-f(-r)-f(-s) \big)
+ \big( f(-r+s)-f(-r)-f(s)  \big)$$
$$
+ \big( f(r-t) -f(r)-f(-s)  \big)
+ \big( f(r+s) -f(r)-f(s)   \big).
$$
Further,
$$
\frac{\partial^2}{\partial r \partial s }\phi(r,s)=
f''(-r-s)+f''(r+s)-f''(-r+s)-f''(r-s)\ge 0
$$
for positive $r,s$ because $f''(t)+f''(-t)$ is nondecreasing for positive $t$.

Thus
$$0\le \int_0^b \int_0^a
\frac{\partial^2}{\partial r \partial s }\phi(r,s)
\,dr\, ds =
\phi(a,b)-\phi(a,0)-\phi(0,b)+\phi(0,0)=\phi(a,b).
$$
Therefore $\phi(a,b)\ge 0$,
and hence the theorem is proved.

\subsection{Proof of Corollary 1}

In the case $1<\rho<2$ the function $f(y)=-|y|^\rho$ satisfies the conditions of Theorem~1.
Hence, (\ref{c11}) is valid for $1<\rho<2$.

It remains to show that the left inequality in (\ref{c11}) holds for $0<\rho<1$. By Theorem~E it suffices to show that, for any $a$ and $b$, $0<a<b$,
$$\phi(a,b):=(a+b)^\rho + (b-a)^\rho - 2^{\rho-1}a^{\rho}-2^{\rho-1}b^{\rho}\ge 0.$$
We have $\phi(a,b)=a^\rho h(z)$, where $z=b/a$,
$$h(z)=(1+z)^\rho + (z-1)^\rho - 2^{\rho-1} - 2^{\rho-1} z^{\rho},$$
Thus $h(z)\ge0$ for $z\ge 1$ because
$h(1)=0$, $h'(z)>0$ for $z\ge 1$.

Hence, (\ref{c11}) is valid for $0<\rho<1$.

The corollary is proved.

\section*{Appendix A. A simple proof of Theorem B}
I.S. Borisov in an oral conversation has proposed the following proof of Theorem B.
We have
$$f(\xi+\eta)=f(\xi)+f'(\xi)\eta+\eta^2\int_0^1(1-z)f''(\xi+z \eta)dz.
$$
Now let us consider the expectation of the last integral.
By convexity of $f''$,
$$\E\left(\eta^2\int_0^1(1-z)f''(\xi+z \eta)dz\ \Big|\ \eta\right)\ge
\eta^2\int_0^1(1-z)f''(z \eta)dz=
f(\eta) - f(0) - f'(o)\eta.
$$
Thus
$$\E f(\xi+\eta)\ge \E f(\xi) + \E f(\eta) + \E f'(\xi)\eta - \E f'(0)\eta.
$$

The restriction of this proof is that all the needed moments, such as $\E \eta^2 f''(\xi+z \eta)$
for $0<z<1$, must exist. For, roughly speaking, ``regular'' functions $f(y)$ growing not faster than $e^{c|y|}$, $c=$\,const, the existence of the moments
follows from
the monotonicity of the functions $f''(y)$,
$f'(y)$, $f(y)$ for sufficiently large $y$ and for sufficiently large~$-y$.

As was noted above, for independent $\xi$ and $\eta$, these moment restrictions can be omitted by virtue of Theorem~E.


\end{document}